\documentclass[10pt,a4paper]{article}

\newcommand{\cH}{\mathcal{H}}

\newcommand{\T}[1]{\mathbb{T}^{#1}}

\newcommand{\R}[1]{\mathbb{R}^{#1}}
\newcommand{\Z}[1]{\mathbb{Z}^{#1}}
\newcommand{\C}[1]{\mathbb{C}^{#1}}

\renewcommand{\H}[1]{\mathbb{H}^{#1}}
\newcommand{\p}{\partial}

\newcommand{\ins}[3]{\vbox to0pt{\kern-#2 \hbox{\kern#1
#3}\vss}\nointerlineskip}

\newcommand{\vf}{\varphi}

\usepackage{graphicx}
\usepackage{bm}
\usepackage{color}
\usepackage{amsmath,amssymb,amsthm,amscd}
\usepackage{diagrams}

\title{Hamiltonian monodromy
via geometric quantization and theta functions}
\author{ Nicola Sansonetto and Mauro Spera \\
Dipartimento di Informatica, Universit\`a degli Studi di Verona\\
Ca' Vignal 2, Strada Le Grazie 15, 37134 Verona\\
{\small e-mail addresses:  nicola.sansonetto@gmail.com, mauro.spera@univr.it}}
\date{}
\begin{document}
\maketitle
\theoremstyle{remark}
\newtheorem{definition}{\bf Definition}
\newtheorem{remark}{\bf Remark}
\newtheorem*{remarks}{\bf Remarks}
\newtheorem{theorem}[definition]{\bf Theorem}
\newtheorem{proposition}[definition]{\bf Proposition}
\newtheorem{corollary}[definition]{\bf Corollary}
\newtheorem{example}{\bf Example}
\newtheorem{exercise}{\bf Exercise}
\newtheorem{application}{\bf Application}
\newtheorem{lemma}[definition]{\bf Lemma}
\newtheorem{fact}{\bf Fact}
\newtheorem{N.B.}[definition]{\bf N.B.}
\null
\null
\centerline{{\bf Abstract}}
In this paper, Hamiltonian monodromy is addressed from  the point of view of
geometric quantization,  and various differential geometric aspects thereof are dealt with, all related to holonomies of 
suitable flat connections.
In the case of completely integrable Hamiltonian systems with two degrees of freedom, 
a link is established between  monodromy and ($2$-level) theta
functions, by resorting to the by now classical differential geometric intepretation of the latter
as covariantly constant sections of a flat connection, via the heat equation.
Furthermore, it is shown that monodromy is tied to the braiding of the Weiestra\ss \ roots
pertaining to  a Lagrangian torus, when endowed with a natural complex structure 
(making it an elliptic curve) manufactured from a natural basis of cycles thereon.
Finally, a new derivation of the monodromy of the spherical pendulum is provided.\par

\section{Introduction}
In this paper,
acting within  the framework of Bohr-Sommerfeld and K\"ahlerian geometric quantization, 
we discuss classical and quantum monodromy from several viewpoints, all related to parallel transport via suitable flat connections. Monodromy, together with the so-called Chern-Duistermaat class, provides an obstruction to the global definition of action-angle variables for completely integrable Hamiltonian systems (\cite{duistermaat, N}; see Subsection 2.1 for details; we do not deal with the non-commutative case, for which we refer to \cite{dazord-delzant, fasso}).  
Our specific contributions consist, first af all, in reinterpreting the Ehresmann-Weinstein connection 
arising from the traditional treatment (see e.g. \cite{duistermaat, bates, cushman-bates}) in vector bundle terms.
Subsequently, we relate monodromy to the freedom of choice of a prequantum connection, and in particular we find that it may be viewed as the obstruction to patching together geometric prequantization bundles equipped with local ``BS-adapted" connections (see Section 3 for precise definitions).  Also, we discuss it
in relation to ${\cal G}_0$-equivalence
of connections (connected component of the identity of the gauge group ${\cal G}$ of a prequantum line bundle), showing, in addition, that it can be detected via
a shift of the quantum action operators (constructed via the recipe of geometric quantization), see Theorem 4. Indeed, in experiments,  monodromy manifests
itself via a shift of the energy levels (\cite{child,CWT}).
Moreover, in the case of completely integrable Hamiltonian systems with two degrees of freedom, we further
relate monodromy  to theta function theory, via the
differential geometric interpretation of the heat equation fulfilled by the $k$-level theta functions
going back to \cite{Witten,adpw,hitchin}.
More precisely, Theorem 5
shows the existence of a representation of the fundamental group $\pi_1(B)$ 
of the base space $B$ of the Lagrangian fibration in tori pertaining 
to a Hamiltonian completely integrable system with two degrees of freedom,
via the holonomy of a flat connection living on a natural complex vector bundle 
(of rank 2) made up of the (2-level) theta functions (pulled back) over $B$. 
The non triviality of this representation signals the emergence of
monodromy.
The upshot is that monodromy can be read via a Berry-type phase shift  
on the space of theta functions of level $2$, manifesting itself as a ``phase gate" (see e.g. \cite{CJ} and Section 4).
The appearance of theta functions in this context is quite natural  from a mechanical
point of view: briefly, this goes as follows.
Given a  basis of cycles on a Liouville torus, constructed as in \cite{Ngoc1,Ngoc2} 
- see also Section 4 below - a natural
{\it complex structure} thereon is determined upon setting
\[
\tau = -\Theta + i T ,
\]
where $\Theta$ is the rotation number and $T > 0$ is the (Poincar\'e) first return time of a point on one of the basis cycles, denoted by $\gamma_1$ (the rotation number is essentially the discrepancy, 
measured on the cycle $\gamma_1$ - corresponding to one of the actions -
between the final and initial position of the aforementioned point, and one can easily manufacture a cycle $\gamma_2$ from such an arrangement;
we notice that in \cite{Ngoc1,Ngoc2}, the roles of the $\gamma$'s are interchanged).
Therefore, each Liouville torus comes equipped with a polarization making it an abelian variety, and hence with a K\"ahler structure (it goes without saying that
the original symplectic form vanishes when restricted to a Liouville torus), and we have a family
of (unobstructed) geometric quantizations of such tori, yielding precisely the theta functions of level $k$
as their quantum Hilbert space (also, they can be adjusted so as to yield orthonormal bases thereof).
Upon varying $\tau$ on the Poincar\'e upper half-plane $\H{}$, one gets a vector bundle whose generic fibre
is given by  the $2$-level theta functions, which has a natural flat connection for which the latter
 are the covariantly constant sections. The ensuing parallel transport translates into
the heat equation fulfilled by the thetas.
This ``universal" construction, pulled back to $B$ via the local action variable map,
yields the above mentioned flat connection, which incorporates monodromy (Theorem 5).\par
 We also point out the direct relationship between the variation of the rotation number
(producing monodromy) and the braiding of the Weierstra\ss \ roots of the elliptic curve associated to $\tau$, again via theta functions; also, a possibly new quick derivation of the monodromy of the spherical pendulum (see e.g. \cite{duistermaat,cushman-bates}), is devised, relying on the above techniques. \par
The present work is organized as follows.
In Section  2 we first collect some background material on monodromy, with special emphasis on the
 two degrees of freedom case - where some 
simplifications occur, notably the vanishing of the Chern-Duistermaat class 
(cf. \cite{duistermaat}) and, what is crucial for our analysis, the existence
of the rotation number (see e.g. \cite{cushman-bates}) - 
and subsequently addressing geometric quantization, focussing our attention on gauge equivalence of connections and reviewing the Bohr-Sommerfeld
conditions, together
 with a brief 
discussion of Hitchin's treatment of polarization independence 
tailored to our purposes (\cite{hitchin}). Also, we give a short account of basic theta function theory and 
its relationship with elliptic curves in Weierstra\ss \ form. The discussion of new results starts in Section 3.
First we discuss the various differential geometric aspects of monodromy hinted at above.
In Section 4 we deal with the theta function approach previously illustrated, and
we establish the relationship between the variation of the rotation number
and the braiding of the Weierstra\ss \ roots of the elliptic curve associated
to $\tau$.  In Section 5 we derive
 the monodromy of the spherical pendulum (see e.g. \cite{duistermaat,cushman-bates}), 
by analysing  suitable elliptic integrals of the first and third kind, and (Section 6) we close the paper with some final remarks and outlook.\par

\section{An overview of integrable systems, geometric quantization and theta functions}
In this section we review some basic facts about completely integrable Hamiltonian 
systems, geometric quantization and theta function theory, for the sake of readability. We will also introduce the notation
that will be used throughout the paper.\par
\subsection{Completely integrable Hamiltonian systems}\label{sec:CIHS}
Let $(M,\omega)$ be a $2n$-dimensional symplectic manifold, and fix
$h:M\longrightarrow \R{}$, a smooth function on $M$ (the Hamiltonian),
with  its associated vector field $X_h$,  fulfilling $i_{X_h}\,\omega = -dh$.
The triple  $(M,\omega, h)$ is called a Hamiltonian system on $M$, with $n$ degrees of freedom,
and it is said to be completely integrable if it 
admits $n$ mutually Poisson-commuting first integrals, which are linearly independent 
almost everywhere in $M$, and, restricting the latter, if necessary, the 
joint level sets of the first integrals are compact and connected.
The Liouville-Arnol'd 
Theorem (see e.g. \cite{arnold,AG,duistermaat}) gives sufficient conditions for the complete 
integrability of a Hamiltonian system.
\begin{theorem}[Liouville-Arnol'd]\label{sec:LA}
 {\it Let $(M,\omega)$ be a $2n$-dimensional symplectic manifold. Let 
  $f = (f_1,\dots,f_n):M\longrightarrow\R{n}$ be a surjective submersion (i.e. the {\rm energy-momentum mapping}),
  such that its components pairwise Poisson-commute. Let $B$ be the set of regular values of $f$. 
  Then for each $b \in B$:}
  \begin{itemize}
    \item[1.] {\it the compact and connected components $f_c^{-1}(b)$ of
                   $f^{-1}(b)$ are diffeomorphic to $\T{n}$;}
    \item[2.] {\it there exists an open neighborhood $U_b$ 
                of $b$ in $B$ and a diffeomorphism 
                \begin{equation}
                   ({\bm I},{\bm \varphi}): f^{-1}(U_b)\longrightarrow V\times\T{n}
                \end{equation}
                with $V$ an open subset of $\R{n}$ such that ${\bm I} = (I_1,\cdots,I_n) = \kappa\circ f$ for
                some diffeomorphism $\kappa :f(U_b)\longrightarrow V$. }
  \item[3.] {\it The coordinates $({\bm I},{\bm \varphi})$ on $M$ are Darboux coordinates, that is 
                 \begin{equation}
                   \omega = d\bm I\wedge d\bm\vf
                 \end{equation}      }           
  \end{itemize}                 
\end{theorem}
\noindent
where ${\bm I}$  is regarded, for future use also, as a row vector, whereas  ${\bm \varphi} = (\varphi_1,...\varphi_n)^T$ is a column vector (see also Subsection 3.2).
From a geometric point of view the Liouville-Arnol'd Theorem 
ensures that $M$ has a $\T{n}$-bundle structure with Lagrangian fibres;
moreover, at the (semi-)local level $f^{-1}(B)$ 
is a Lagrangian toric principal bundle with structure group $\T{n}$, 
the fibres are Lagrangian and the structure group acts in a Hamiltonian way, 
with momentum map given by the projection bundle map. The construction of the 
toric principal bundle or, equivalently, the existence
of global action-angle coordinates is only (semi-)local; indeed,
Duistermaat proved the following:
\begin{theorem} (\cite{duistermaat})
  {\it The $\T{n}$-bundle $\pi:f^{-1}(B)\longrightarrow \R{n}$ 
  is topologically trivial if and only if the monodromy and the Chern-Duistermaat class 
  of the $\T{n}$-bundle are trivial.
  Moreover if the symplectic form is exact then the existence of global 
  action-angle coordinates is equivalent to the triviality of the Lagrangian toric fibration.}
\end{theorem}  
See also \cite{N}.
\begin{remarks}
  \begin{enumerate}
    \item Geometrically, monodromy is the obstruction preventing
              the $\T{n}$-bundle from being a principal bundle with structure group $\T{n}$,
              whilst the Chern-Duistermaat class is the obstruction to the existence of a global section
              of the $\T{n}$-bundle.
    \item Observe that in the case of a system with two degrees of freedom  possessing an isolated critical 
              value (of focus--focus type) of the energy-momentum map $f$, the Chern-Duistermaat class is trivial since
              $B$ admits a Leray cover with empty triple intersections. Therefore the only
              obstruction to the triviality of the fibration is monodromy. We shall assume this condition in the sequel.
  \end{enumerate}              
\end{remarks}

Zung (\cite{Z}) gives a sufficient condition for the non-triviality of monodromy near isolated focus-focus singularities:
more precisely, the (local) monodromy near a topologically stable focus-focus point (in the interior of the energy-momentum
range) is non-trivial. This result will be used in Section 4.\par

It will be convenient for us to study Hamiltonian
monodromy from a differential geometric point of view 
(see \cite{weinstein,duistermaat,cushman-bates}). Indeed it is well-known (\cite{weinstein})
that a Lagrangian fibration admits an affine, flat, torsion free connection
$\nabla^{Ehr}: TM\longrightarrow VM$ -  the vertical bundle over $M$ -
on the Lagrangian leaves, which is an Ehresmann good connection for the fibration (i.e. that is every smooth
curve on the base has a horizontal lift). 
The  $GL(n, \Z{})$-holonomy representation  $hol(\nabla^{Ehr})$
of $\nabla^{Ehr}$ is the monodromy representation 
 of the $\T{n}$-bundle $\pi: f^{-1}(B)\longrightarrow  \R{n}$, therefore if the monodromy 
is non trivial, then the $\T{n}$-bundle is not principal. Moreover the monodromy representation actually
takes values in  $SL(n,\Z{})$ upon choosing  suitable bases of the tangent spaces
of the base space. In Subsection 3.1 we will reformulate the above discussion in vector bundle terms.

\subsection{Geometric quantization}\label{sec:GQ}
Let us now briefly review the basics of geometric quantization; we refer to \cite{woodhouse,Bry,kirillov,souriau,kostant}
for a complete account.
Recall that if $(M,\omega)$ is a symplectic manifold of (real) dimension $2n$ 
such that $\left[\frac{1}{2\pi}\omega\right]\in H^2(M,\Z{})$, then the Weil-Kostant Theorem  
states that there exists a complex line bundle $(L,\nabla,h)$ over $M$ equipped 
with a hermitian metric $h$ and a compatible connection $\nabla$ with curvature 
$F_{\nabla} = \omega$.
Hence $[\omega] = c_1(L)$, the first Chern class of $L\rightarrow M$. The connection 
$\nabla$ is called a prequantum connection and $L\rightarrow M$ the 
prequantum line bundle. The different choices of $L\rightarrow M$ and $\nabla$ are parametrized
by the first cohomology group $H^1(M,S^1)$ (see e.g. \cite{woodhouse}, Ch.8).
In more detail (see  also \cite{Mo}, 1.7), given any complex line bundle $L\rightarrow M$,
the connections thereon are classified, up to gauge equivalence, by their
{\it curvature} (fixing the topological type of the line bundle, via the first Chern class) and
by their {\it holonomy}, specified, in turn, on a basis of (real) homology 1-cycles $[\gamma_i]$,
for $H_1(M, \R{})$, of dimension $b_1$, the first Betti number of $M$ - represented, for instance,
by smooth curves passing through a given point. The holonomy is trivial if  $M$ is simply
connected.  The gauge group ${\cal G}$ consists, in this case, of
all smooth maps $g: M \rightarrow S^1$ - explicitly, $ g: x \mapsto \exp [{i \, \varphi}(x)] $,
obvious notation - and it
is not connected in general, its connected components being parametrized 
by the degree of the maps $g: M \rightarrow S^1$.
The connected component (of the identity) of ${\cal G}$ will be denoted by ${\cal G}_0$, as usual, and will
play an important role in what follows.\par
\par
Given a connection $\nabla_0$, any other connection is of the form 
$\nabla = \nabla_0 + a$, with $a \in \Lambda^{1}(M)$, 
(i.e. they build up an affine space modelled on the space of 1-forms $ \Lambda^{1}(M)$) and the relation between their respective
curvatures is
\begin{equation}
F_{\nabla } = F_{\nabla_0 } + da
\end{equation}
Therefore, the curvatures are the same if and only if $a$ is closed. This being the case, $a$ determines a
de Rham cohomology class $[a] \in H^{1}(M,\R{})$,  fully recovered via the {\it period map}
\begin{equation}
H^{1}(M,\R{}) \ni [a] \mapsto \left(\int_{\gamma_1} a\, , \dots , \int_{\gamma_{b_1}}  a \,\right) \in {\R{}}^n
\end{equation}
The gauge group ${\cal G}$ acts on connections via
\begin{equation}
\nabla \mapsto \nabla + g \cdot d \, g^{-1} = \nabla  - i \, d \varphi 
\end{equation}
Therefore, the set of all gauge inequivalent connections (possessing the same curvature) is clearly given by
\begin{equation}
H^{1}(M,\R{})/ H^{1}(M,\Z{})
\end{equation}
and, if $M$ is a {\it torus}, then the above set is again a torus, the {\it Jacobian} of $M$.
If the initial connection has zero curvature, then the above space parametrises 
{\it flat} connections up to gauge equivalence.\par
Coming back to the specific geometric quantization setting,
given a Lagrangian submanifold $\Lambda$ of the symplectic manifold $M$, the symplectic 2-form $\omega$
vanishes upon restriction to $\Lambda$ by definition, and any (semi-local) symplectic
potential $\theta$ becomes a closed form thereon, defining a (semi-local) connection form
pertaining to the restriction of the prequantum connection $\nabla$, denoted by the
same symbol. The latter is a {\it flat} connection and a global covariantly  constant 
section of the restriction of the prequantum line bundle exists if and only if it has 
trivial holonomy, that is, the induced character $\chi:\pi_1(\Lambda)\longrightarrow U(1)$
is trivial (see e.g. \cite{tyurin1}), or, equivalently, that the {\it Bohr-Sommerfeld} condition is fulfilled:
\begin{equation}
\left[\frac{1}{2\pi}\theta\right]\in H^1(M,\Z{}) \quad\hbox{i.e.}\quad \int_\gamma \theta\, \in\, 2\pi\Z{}
\end{equation}
for any closed loop $\gamma$ in $\Lambda$. 
\par
 A covariantly constant section (which we call WKB-, or BS-wave function) takes the form
\begin{equation}\label{eq:WKB}
  s(m) := hol_\gamma(\nabla)\,\cdot\, s(m_0) = e^{i\int_\gamma \theta}\, s(m_0)
\end{equation}
with $\gamma$ denoting any path connecting a chosen point $m_0$ in
$\Lambda$ with a generic point $m\in \Lambda$, $hol_\gamma(\nabla)$ being
the holonomy along $\gamma$ of the restriction to $\Lambda$ of the prequantum connection $\nabla$.
The r.h.s. of (\ref{eq:WKB}) tacitly assumes the choice of a trivialization 
of $L\mid_{\Lambda}\longrightarrow \Lambda$ around $m_0$ and $m$ in a corresponding local chart.

\begin{remarks}
  \begin{enumerate}
    \item We stress the fact that the Bohr-Sommerfeld condition forces us to deal with ${\cal G}_0$-equivalence
    classes (i.e. the degree of the gauge maps must be zero) in order to avoid trivialities. See in particular Subsection 3.3.
    \item Our definition of WKB-wave function is slightly different from the conventional one
              (see e.g. \cite{woodhouse}). Indeed we do not require square-integrability and we 
              do not twist the prequantization bundle with $\Delta_\nabla$ (whose smooth sections
              consist of the complex n-forms on $\Lambda$), thus neglecting the 
              ``amplitude-squared''.
             \item There is a version of the Bohr-Sommerfeld condition incorporating the Maslov class,
               but we shall not need this refinement in what follows. 
  \end{enumerate}
\end{remarks}
We shall resume the above discussion in Section 3.\par
We also recall that the prequantum  connection $\nabla$ allows the construction of the 
(Hermitian) prequantum observables $Q(\cdot)$
via the formula
\begin{equation}
 Q(f) = - i\nabla_{X_f}  +  f = -i X_f - i_{X_f} \theta  +  f  
\end{equation}
The connection is determined up to a closed 1-form,  yielding a corresponding ambiguity
in the definition of the quantum observable $Q(f)$ attached to $f$. This fact will be exploited in the sequel
(see again Section 3).\par

In the K\"ahler case one can perform  {\it holomorphic quantization}, whereby one
takes the space of holomorphic sections $H^{0}(L, J)$ of a holomorphic prequantum line bundle,
provided it is not trivial, as the Hilbert space of the theory
($J$ denotes a complex structure on $M$, see e.g. \cite{hitchin} for details). In this case there is a 
canonically defined connection, called the Chern, or Chern-Bott connection, compatible 
with both the hermitian and the holomorphic structure (cf. \cite{GH}). Independence of polarization
(i.e. of the complex structure, in this case) is achieved once one finds a (projectively)
flat connection on the vector bundle $V \rightarrow {\cal T}$ with fibre $H^{0}(L, J)$
(of constant dimension, under suitable assumptions provided by the Kodaira vanishing theorem)
over the (Teichm\"uller) space of complex structures ${\cal T}$.
An important example, which will be
needed later on, is provided by the $k$-level theta functions, which can be viewed as (a basis of) the space of holomorphic sections of a holomorphic line bundle
(the $k$th tensor product of the theta line bundle) defined on a principally polarized abelian variety  (\cite{Witten, adpw,hitchin}, see also \cite{Spe}). 
It follows from the Riemann-Roch theorem that this space has (complex) dimension $k$.
In dimension two, the role of ${\cal T}$ is played by
the Poincar\'e upper half plane $\H{}$ (a complex structure being labelled by $\tau \in \H{}$).
The covariant constancy of the thetas is ascribed to their fulfilment of the heat equation.
In the following section we give some extra details on theta functions needed for the sequel.

\subsection{Elliptic integrals and theta functions}\label{sec:theta}
In this subsection we collect some facts about elliptic integrals and theta functions in one variable, in view of future use. The theory is thoroughly expounded in
many classical texts, see e.g. among others  \cite{Mum,McK-Moll,Tricomi,WW,
GH,Kempf}. We shall use this material in Sections 4 and 5.\par
\smallskip
The Weierstra\ss \ canonical forms of the elliptic integrals of the first, second and third kind read, respectively:\par

\begin{equation}
I_1 = \int {dz \over {\sqrt{P(z)} }}, \quad
I_2 = \int {z\,dz \over {\sqrt{P(z)} }}, \quad
I_3 = \int {dz \over {(z - c) \sqrt{P(z)} }}
\end{equation}
where
\begin{equation}
 P(x) \, := \, 4x^3 - g_2 x - g_3 \, = \, 4(x - e_1)(x - e_2)(x - e_3)
\end{equation}
with $ e_1 + e_2 + e_3 = 0$ (the $e_i$'s are all distinct); in $I_3$, $c$ is required not to be a root of $P$.
 The elliptic integral $I_1$ above is explicitly inverted
by the celebrated Weiestra\ss  \ function $\wp = \wp (z, g_2, g_3) \equiv \wp (z, \tau)$,
fulfilling $ y^2 = P(x) $, giving rise to an elliptic curve  ${\cal C}$, with $x = \wp$, $y = {\wp}\prime$.
Then ${\cal C} \cong
 \C{}/ \, \Z{} + \Z{}\tau$, the  torus defined by quotienting $\C{}$ by a normalized lattice $\Z{} + \Z{}\tau$,
 where $\tau =  {{\omega^{\prime} }\over {\omega}} \in \C{}, \Im{\tau} > 0)$ (ratio of (half)-periods).
One has $ e_i = \wp(\omega_i)$, where  $\omega_1 = \omega$,  $\omega_2 = \omega + \omega^{\prime}$,  $\omega_3 = \omega^{\prime}$.
The (Jacobi) modulus (squared) $k^2$  (with $k \in \C{}\setminus \{  0,1 \}$) together with its complementary modulus $k^{\prime}$ fulfilling ${k^{\prime}}^2 = 1 - k^2$,  
can be interpreted as the simple ratio of the three roots of  $P$ (see below).
The standard theta function reads
\begin{equation}
\vartheta (z, \tau ) = \sum_{n \in \Z{}} e^{i\pi\, n^2 \tau  + 2\pi i \, n z  }
\end{equation}
Let us also  record the expressions for theta function with $2$-characteristics  (using Mumford's notation (\cite{Mum}):
\begin{equation}
\vartheta _{ab}(z, \tau ) =    e^{\pi i a^2 \tau + 2\pi i a (z + b)}\vartheta( z + a \tau + b, \tau)
\end{equation}
where $a, b \in {1\over 2} \Z{}$. Comparison with traditional notations yields
$\vartheta _{00} = \vartheta _{3}$, $\vartheta _{0{1\over 2}}  \equiv \vartheta _{01} = \vartheta _{4}$, 
$\vartheta _{{1\over 2}0}  \equiv \vartheta _{10} = \vartheta _{2}$, 
$\vartheta _{{1\over 2}{1\over 2}}  \equiv \vartheta _{11} = \vartheta _{1}$

The Jacobi modulus $k$ of the attached elliptic curve can be recovered from $\tau$ via the formula
\begin{equation}
k^2 = {{{\vartheta_2}^4 (0, \tau)}\over {{\vartheta_3}^4 (0, \tau)} } = {{e_2 - e_3} \over {e_1 - e_3} }
\end{equation}
(this is the very motivation which led Jacobi to devising theta functions).
\par

Indeed, let us recall, for future use, the following expressions relating the Weierstra\ss \ roots to theta functions:

\begin{equation}
e_2 - e_3 = \left( {{\pi}\over {2\omega}} \right)^{{}^2}{\vartheta_2}^4 (0, \tau)
\qquad
e_1 - e_2 = \left( {{\pi}\over {2\omega}} \right)^{{}^2}{\vartheta_4}^4 (0, \tau)
\end{equation}
following from
\begin{equation}
\begin{matrix} 
e_1 & =  & {{\pi}^2\over {12\,\omega^2}}[ {\vartheta_3}^4 (0, \tau) + {\vartheta_4}^4 (0, \tau)] \cr

e_2 & =  &{{\pi}^2\over {12\,\omega^2}}[ {\vartheta_2}^4 (0, \tau) - {\vartheta_4}^4 (0, \tau)]  \cr

e_3 & = & - \, {{\pi}^2\over {12\,\omega^2}}[ {\vartheta_2}^4 (0, \tau) + {\vartheta_3}^4 (0, \tau)], \cr
\end{matrix} 
\end{equation}
the beautiful Jacobi formula:
\begin{equation}
{\vartheta_2}^4 (0, \tau) + {\vartheta_4}^4 (0, \tau) = {\vartheta_3}^4 (0, \tau)
\end{equation}
and, most important, the following transformation law:
\begin{equation}
\begin{matrix}
\vartheta_1 (z, \tau + 1) & = & e^{i {\pi \over 4}} \, \vartheta_1 (z, \tau) \cr
\vartheta_2 (z, \tau + 1) &= &e^{i {\pi \over 4}} \, \vartheta_2 (z, \tau)  \cr
\vartheta_3 (z, \tau + 1) & = &\vartheta_4(z, \tau) \cr
\vartheta_4 (z, \tau + 1) & = &\vartheta_3(z, \tau) \cr
\end{matrix}
\end{equation}

Let us now consider the following modified theta function:
\begin{equation}
{\widetilde\vartheta} (z, \tau )  = e^{{\pi\over 2}\, ({\Im \tau})^{-1}{z^2}} \vartheta (z, \tau )
\end{equation}
Notice that the prefactor $  e^{{\pi\over 2}\, ({\Im \tau})^{-1}{z^2}} $ is invariant with respect to the
trasformation $\tau \mapsto \tau +1$. It is this modified theta function that, in the algebro-geometric literature
(see e.g. \cite{Kempf}) gives rise to the
unique (up to a constant) 
holomorphic section of the theta line bundle associated to a complex torus (and, in general, to a principally polarized abelian variety), which is 
actually the prequantum bundle  (\cite{Witten, adpw,hitchin}, see also \cite{Spe}). This is readily generalized to the $k$-level theta functions, which
(up to constants) yield an orthonormal basis for the ($k$-dimensional, by Riemann-Roch) quantum Hilbert space (see  \cite{Kempf, Loi}).\par
We record the relevant formulae, for definiteness (with a slightly different notation, also in order to avoid confusion with theta functions with characteristics):
\begin{equation}
\widetilde{\theta}_{k,j} (z, \tau) = e^{k{\pi\over 2}\, ({\Im \tau})^{-1}{z^2}} \sum_{n \in \Z{}} e^{{i \over k}\pi\, (k n + j)^2 \tau  + 2\pi i \, (k n + j) z  } \equiv 
e^{k{\pi\over 2}\, ({\Im \tau})^{-1}{z^2}} {\theta}_{k,j} (z, \tau) 
\end{equation}
for $j = 0,\dots , k-1$.
A crucial fact for what follows is that the $k$-level theta functions ${\theta}_{k,j}$ fulfil the (holomorphic) {\it heat equation}
\begin{equation}
\left[{{\partial} \over {\partial\tau}} + {1\over{4\pi k}}\, {{\partial^2} \over {\partial z^2}} \right]\,{\theta}_{k,j}  = 0
\end{equation}
Now, a straightforward computation shows that, under the trasformation $\tau \mapsto \tau + 1$, the
$2$-level theta functions $\vartheta_{2,0} $ and $\vartheta_{2,1} $, together with their ``tilded"  analogues, behave as follows
\begin{equation}
\theta_{2,0} (z, \tau + 1) =  \theta_{2,0} (z, \tau) , \qquad  \theta_{2,1} (z, \tau + 1) =  e^{i{\pi\over 2}}\, \theta_{2,1} (z, \tau) 
\end{equation}
Consider the vector bundle $V \rightarrow \H{}$, with $V_{\tau}$ 
(fibre at $\tau$) given by the 2-dimensional complex vector space of  $2$-level theta functions with fixed parameter $\tau$. It comes equipped
  with the {\it heat connection} $\nabla$, and the $2$-level  theta functions provide
 a basis of covariantly constant sections thereof, this being expressed by fulfilment of the
 heat equation. An important consequence is that, in particular,
 the natural $SL(2, {\bf \Z{}})$-action on $\H{}$ given by
 \begin{equation}
 \tau \mapsto {{a\, \tau +  b}  \over {c\, \tau +  d}}
 \end{equation}
 ($ad - bc = 1$),
 yields, in turn,
a parallel displacement map $Q(Z): V_{\tau} \rightarrow V_{Z \cdot\tau}$, for $Z \in SL(2, \Z{})$ (along any path connecting the two points). Specifically, for the 
matrix $Z_0$ associated to the map $\tau \mapsto \tau + 1$, i.e.
 \begin{equation}
Z_0 = \begin{pmatrix} 
         1 & 1  \cr
         0 & 1  \cr
        \end{pmatrix}  
 \end{equation}
one has the (``phase gate" \cite{CJ}) matrix, whereby we rephrase the transformation formula for the $\theta_{2,j}$'s  and the $\widetilde{\theta}_{2,j}$'s :
 \begin{equation}
Q(Z_0)  \, = \, \begin{pmatrix} 1 & 0   \cr
         0 &  e^{i{\pi\over 2}}\cr
         \end{pmatrix}        
\end{equation}
acting on the theta vector
 $ {\underline{\widetilde{\theta}}}_2 (z) = ( \widetilde{\theta}_{2,0} (z, \tau), \widetilde{\theta}_{2,1} (z,\tau))^T$. 
 Notice that $Q(Z_0)^4 = Id_2$.\par

 We also remark that, by virtue of the preceding formulae, the map $\tau \mapsto \tau + 1$ determines a switch of the roots $e_2$ and $e_3$.
 This will be important in Subsection 4.2. \par
 Finally we notice that,
for $\vartheta_3$, one has
 \begin{equation}
\vartheta_3 (0, \tau + 1) = \vartheta_4 (0, \tau) = {{\vartheta_4 (0, \tau )}\over  {\vartheta_3 (0, \tau )}}
\vartheta_3 (0, \tau) = \sqrt{k^{\prime}} \,\vartheta_3 (0, \tau)
 \end{equation}
(by the Jacobi formula), yielding a differential geometric interpretation of the Jacobi modulus.\par

\section{Hamiltonian monodromy and Geometric Quantization}

\subsection{The Weinstein connection revisited}\label{sec:weinstein}
In this Subsection we elaborate, in view of future use, on the canonical connection 
attached to the natural Lagrangian fibration in Liouville tori (\cite{weinstein}),  
see Subsection~\ref{sec:GQ} as well, by rephrasing it in terms of vector bundles.\par
The local action variables ${\bf I } = (I_i)$ provide a local diffeomorphism
between the set of regular values of the moment map $B$ and
 ${\R{}}^n$. On the trivial ${\R{}}^n$-bundle on the latter
one has the natural flat connection induced by $d$. This is
pulled back to $B$, and the local pieces glue together to 
yield a flat connection on the trivial ${\R{}}^n$-bundle
thereon, which we call canonical and denote by ${\nabla}^{can}$, 
whose holonomy $hol(\mathbb{\nabla}^{can})$ (with values in $SL(n, \Z{})$)
is exactly the monodromy (also cf. \cite{duistermaat,bates,cushman-bates,weinstein}). 
The upshot is the following: \par
\begin{proposition}
\begin{enumerate}
\item {\it The local
action variables build up, collectively, a {\rm covariantly
constant section} of ${\nabla}^{can}$.}\par
\item {\it The following relation between {\rm classical} and {\rm quantum} monodromy holds (cf. \cite{Ngoc1,Ngoc2} })
\begin{equation}
\mu_q = (\mu_c)^{-T} 
\end{equation}
\end{enumerate}
\end{proposition}
 Part 2 immediately follows from de Rham's theorem via
\begin{equation}
\int_{\gamma} \theta = \langle [\theta], [\gamma] \rangle = \langle [\theta], Z^{-1}\cdot Z [\gamma] \rangle
= \langle Z^{-T} [\theta], Z [\gamma] \rangle
\end{equation}
(duality pairing between $H_1(\Lambda, \R{})$ and $H^1(\Lambda,\R{})$ and 
 $Z \in SL(n, \Z{})$, via diffeomorphism invariance
$\int_{\gamma} \theta = \int_{\varphi \cdot\gamma } \varphi^{*}\theta$ (obvious notation).\par

\subsection{Monodromy and prequantum connections}
Here we resume the discussion about the freedom of choice of the prequantum connection, by
focussing on the case of Lagrangian fibrations in Liouville tori (cf. Subsection~\ref{sec:CIHS}).
For a trivial Lagrangian bundle  $ U \times {\T{n}} \rightarrow U  $ (actually, its total space),
consider its prequantum line bundle $L \rightarrow U \times {\T{}}^n$, with a prequantizing connection
$\nabla$, with local connection form given by a symplectic potential $\theta$ determined, in a first
instance, up to a closed form. We have two natural choices for the prequantum connection.\par
 Firstly, set
\begin{equation}
\nabla  \leftrightarrow \theta = \sum_{k=1}^{n} I_k d\varphi_k \equiv {\bm I}\, d {\bm \varphi}.
\end{equation}
This may be called BS-adapted (or vertical) connection, since it just comes from a geometrical reformulation of the standard
procedure. It fulfils 
\begin{equation}
\nabla_{X_b} = X_b
\end{equation}
with $X_b$ any vector field on $M$ tangent to a Lagrangian section, and
it is flat along fibres.
More intrinsically,  given an adapted connection as above, the action variables may be recovered as follows:
\begin{equation}
I_k = {1\over{2\pi i}}\log hol(\nabla\!\!\mid_{\Lambda}, \gamma_k)
\end{equation}
where  the $\gamma_k$'s yield a basis of  1-cycles in $\Lambda$, this making their {\it local} character clear.
 Hence, monodromy may be viewed as the obstruction to patching together geometric prequantization bundles equipped with local BS-adapted connections. Of course there is no global obstruction to prequantization {\it tout court}, by Weil-Kostant.\par
Secondly, set
\begin{equation}
\nabla^{\prime}  \leftrightarrow \theta^{\prime} =  -\sum_{k=1}^{n} \varphi_k dI_k
\equiv  - d {\bm I} \, {\bm \varphi} 
\end{equation}
This connection can be termed {\it monodromy connection}, since parallel
transport along a non trivial loop contained in a local Lagrangian section ($\varphi = c$)
(whereupon it is flat)
produces a holonomy given by
\begin{equation}
e^{- c i  \Delta I}
\end{equation}
(obvious notation) tied to the possible non globality of the action variables.
It can be characterised intrinsically as well by the requirement
\begin{equation}
\nabla_{X_{\varphi}} = X_{\varphi}
\end{equation}
(with $X_{\varphi}$ tangent to the fibres).
So the  freedom in choosing the prequantum connection leads to detection of monodromy.\\

 Notice that in the case of $B$ is a multi-punctured domain, its fundamental group
is a free group on $m$ generators (if we have $m$ punctures). Now, the monodromy around
a puncture can be ``signed", so the monodromy representation of suitable non trivial loops
may nevertheless be trivial.

\subsection{Gauge equivalence of flat connections and monodromy}
In this Subsection we further specialise the general discussion outlined in Subsection~\ref{sec:theta} 
and we address monodromy from a gauge theoretic point of view - encompassing
Ngoc's treatment \cite{Ngoc1,Ngoc2}. Let $\Lambda \cong \T{n}$ be a Liouville torus. Then, we have
already noticed, in Subsection~\ref{sec:weinstein}, that, on the one hand,
the homology group $H_1(\Lambda, \Z{})$ is the arena of classical monodromy, stemming
from an $SL(n, \Z{})$-action on the classical cycles. The cohomology group
$H^1(\Lambda, \Z{})$ is, on the other hand, a receptacle for quantum monodromy. \par
Now, what is crucial in highlighting monodromy is  that the finer
notion of ${\cal G}_0$-(in)equivalence should be
used instead of mere gauge equivalence. This goes as follows. Resuming the discussion of 
Subsection~\ref{sec:GQ}, let us take the (integral, upon enforcing BS)
de Rham class of $\nabla \equiv {\nabla}\!\!\mid_{\Lambda}$, i.e. $ [\theta] $, mapping to a point in
${\Z{}}^n$ via the period map
\begin{equation}
H^1(\Lambda, \Z{}) \ni [\theta] \mapsto \left({1\over {2\pi}}\int_{\gamma_1} \theta, ... , {1\over {2\pi}}\int_{\gamma_n} \theta \right)  \in {\Z{}}^n
\end{equation}
(where $(\gamma_1, \dots , \gamma_n)$ is a basis of 1-cycles)
and denote as ${\cal BS}$, for convenience, the set of all classes $[\theta_{\nabla}]$ (it is enough to consider  BS-adapted connections). Then (obvious notation)
\begin{equation}\label{eq:gauge}
{\cal BS} \cong  H^{1}({\T{}}^n, \Z{}) =  {\cal G} \cdot [\nabla_0]
\end{equation}
(with $\nabla_0$ a fixed flat connection). Thus ${\cal BS}$ is  a ${\cal G}$-homogeneous
space $\cong {\Z{}}^n$, whereupon
the group ${\cal G}_0$ acts trivially. Hence
${\cal G}/{\cal G}_0 \cong SL(n,\Z{})$ acts freely on ${\cal BS}$,  and provides the receptacle of  a  natural monodromy representation
\begin{equation}
\widetilde{M}:  \, \pi_1(B) \rightarrow {\cal G}/{\cal G}_0 \cong SL(n, \Z{})
\end{equation}

The BS-wave functions  are, in turn, characters of ${\T{}}^n$, i.e. elements of its dual group, and
the latter is isomorphic to   ${\Z{}}^n$.
Explicitly
one has a family of (flat) BS-connections $\nabla_{ {\bm n}}$, $( {\bm n} \in {\Z{}}^n)$,
(which are all ${\cal G}$-equivalent but not ${\cal G}_0$-equivalent)
with covariantly constant section (up to a constant)
\begin{equation}
s = \chi_{{\bm n} } ({\bm{ \varphi}}) = e^{i  {\bm n} \cdot   {\bm \varphi} }
\end{equation}
We set ${\cal H}_{{\bm n} } = <\chi_{{\bm n} }>$. 
Alternatively, we may proceed as follows and, in order 
to fix ideas, we take $n=2$ and  consider the shift induced by $n_2 \mapsto n_2 + 1$;
it can be ascertained via the following procedure.
Let $-i\partial_{\varphi_2}$ be the quantum observable on ${\cal H}_{{\bm n} } $ associated to $I_2$ (cf. also the more general discussion below) acting via
$-i\partial_{\varphi_2}\chi_{{\bm n} } =  n_2 \chi_{{\bm n} } $.
If $U: {\cal H}_{{\bm n} } \rightarrow {\cal H}_{{\bm m} } $,
with $(m_1, m_2) = (n_1, n_2 + 1)$ is the unitary operator sending $\chi_{{\bm n} }$ to
$\chi_{{\bm m} }$, one finds, on ${\cal H}_{{\bm m} }$,  the shifted operator
\begin{equation}\label{eq:unitaria}
U \circ (-i\partial_{\varphi_2} ) \circ U^{-1} = -i\partial_{\varphi_2}  +  I d
\end{equation}
and this is again a flashing light for monodromy.
If the basis of cycles on a model torus is kept fixed, monodromy can be detected as
a switch to another BS-class: one has a shift
of the action variables
(cf. the parallel transport of lattices defined in \cite{Ngoc1,Ngoc2}).\par

More generally, let us perform a coordinate transformation on a fixed BS-torus 
(under our assumptions we may
neglect the translational part), and let us extend it to a canonical transformation in 
the ambient manifold $M$ (in a fibre neighborhood of the torus in question):
\begin{equation}
{\bm I}^{\prime} = {\bm I}\,  Z^{-1}, \qquad \varphi^{\prime} =  Z \, {\bm \varphi}
\end{equation}
($Z \in SL(n,\Z{})$). One has, indeed
\begin{equation}
d {\bm I}^{\prime} = d{\bm I} \, Z^{-1},\quad  d {\bm \varphi}^{\prime} =  Z \, d  {\bm \varphi},
\quad \partial_{ {\bm \varphi}^{\prime}} = Z^{-T} \partial_{\bm \varphi}
\end{equation}
and
\begin{equation}
 d{\bm I}^{\prime} \wedge d {\bm \varphi}^{\prime} = d {\bm I} \, Z^{-1} \wedge Z \, d{\bm \varphi} =  d {\bm I} \wedge d {\bm \varphi} = \omega 
\end{equation}
Now, the quantum operator  associated to the action variable $I$ according to
the general formula given above is $-i \partial_{\varphi}$ (the last two terms cancel out).
This is checked immediately (obvious notation)
\begin{equation}
\hat{\bm I} =  -i X_{\bm I}  -  i_{X_{\bm I}}\theta +  {\bm I} 
\end{equation}
but $\theta = {\bm I} \, d {\bm \varphi}$ and $X_{\bm I}  =  \partial_{\bm \varphi}$, hence
\begin{equation}\label{eq:osservabile}
\hat{\bm I} =  - i X_{\bm I} \quad \mapsto \quad
\hat{{\bm I}^{\prime}} = -i Z^{-T}\partial_{\bm \varphi}
\end{equation}
and, acting on ${\cal H}_{{\bm n} }$, reproduces the specific result above.
The link with monodromy manifests itself via a non trivial $SL(n,\Z{})$-representation of $\pi_1(B)$
 given by $[\gamma]  \mapsto Z = Z([\gamma])$. It can be viewed as a product $Z = \prod_i Z_i$ of transformations involving two intersecting open charts whereupon no singularity is present. Upon tracing a circuit $\gamma$ surrounding
an isolated singularity (of focus - focus type),
one ends up with the $Z$ above (cf. \cite{N,cushman-bates}).
Summing up, the shift occurs if and only if there is monodromy, and
everything is stored in the commutative diagram below:\par
\begin{equation}
\begin{diagram}
{\cal U} & \rTo^Q & {\cH} & \rTo^{\p_{\bm\vf}} & {\cH}\\
\dTo_{\cal C} & & \dTo_{U} && \dTo_{U}\\
{\cal U} & \rTo^Q & {\cH^{\prime}} & \rTo^{\p_{\bm\vf^{\prime}}} & {\cH}^{\prime} \\
\end{diagram}
\end{equation}

\null
\null
\noindent
where ${\cal U}$ is a fibre neighbourhood of a BS torus, $Q$ is the quantization map producing
the quantum Hilbert space, and ${\cal C}$ is the canonical transformation of the fibre neighbourhood
onto itself described above and $U$ the unitary map connecting  primed and unprimed spaces.\par
Notice that the general formula for the prequantum operator, when applied to ${\bm I}$, should
be appropriately restricted to a fibre neighbourhood of the torus under consideration.
The monodromy action, via $Z$, changes the quantum operator.
Indeed, in  spectroscopy, monodromy manifests itself precisely through 
a shift of the energy levels, see e.g. \cite{child,CWT,EJS} and references therein.\par
\smallskip
The above discussion can be summarised  by means of the following
\begin{theorem} (Gauge theoretic interpretation of monodromy)\par
\begin{enumerate}
\item {\it The monodromy representation  can be viewed as a map
\begin{equation}
\widetilde{M}: \, \pi_1(B) \rightarrow {\cal G}/{\cal G}_0 \cong SL(n, \Z{}),
\end{equation}
which acts transitively on ${\cal BS}$, as expression (\ref{eq:gauge}) shows, 
and can be read both on wave functions and observables.}\par
\item {\it Explicitly, upon
choosing a BS-adapted connection, one can work either with a
fixed basis of cycles,  and then monodromy  induces
a change of connection and Hilbert space in a different 
${\cal G}_0$-class - this however can still be read on a single Hilbert space, cf. (\ref{eq:unitaria})) - or, alternatively,  with a change of coordinates, remaining 
in the same Hilbert space,  causing eventually 
a change in the quantum action operator, (\ref{eq:osservabile}).}
\end{enumerate}
\end{theorem}

\section{Hamiltonian monodromy via theta functions}

 From now on we confine ourselves to completely integrable Hamiltonian systems with two
degrees of freedom.\par
\subsection{The heat connection}
As we have already seen in Section 2, the 2-torus bundle $f: M\longrightarrow B$
has monodromy if and only is the holonomy of the Ehresmann connection 
$\nabla^{Ehr}$ on $B$ is non-trivial. Now we can relate the monodromy 
of the fibration $f$ to the holonomy of the {\it heat connection} introduced right below.\par
Define a map $\tau_{U} : B \supset U \rightarrow \H{}$ via
 $\tau (b) := -\Theta (b) + i \, T (b) $ (notice that $\Im (\tau) > 0$)
 using a basis $(\gamma_1, \gamma_2)$ for the cycles as in \cite{Ngoc1} (with the roles of the $\gamma_i$ interchanged, also cf. Introduction). Note that this is the crucial point wherein two-dimensionality intervenes.\par
Resuming the 2-level theta vector bundle $V \rightarrow \H{}$,  one constructs the pulled-back bundle
 \begin{equation}
\tau_{U}^{*}V \rightarrow U
 \end{equation}
equipped with a flat connection $\nabla_U = \tau_{U}^{*}\nabla$ ($\nabla$ is the ``old" heat connection on the theta bundle). Gluing these local bundles together
one ends up with
a (smooth) vector bundle ${\cal V} \rightarrow B$,
again endowed with a flat connection, called again heat connection and denoted by $\nabla^{heat}$. Clearly, in view of the discussion
in Subsection 2.3, the following holds 
\smallskip
\begin{theorem}
  {\it  Let $(M,\omega,h)$ be a completely integrable Hamiltonian system with two degrees of freedom, possessing a finite number of singularities of focus-focus type
  (cf. Subsection 2.1).
  Then}
\begin{enumerate}
 \item {\it The holonomy of the heat connection on ${\cal V} \rightarrow B$, the pulled-back $2$-level theta vector bundle, relates
  to the holonomy of the canonical connection (Subsection 3.1) in the following guise}
  \begin{equation} 
    hol(\nabla^{heat}) \, = Q (hol(\mathbb{\nabla}^{can}))
 \end{equation}
   \item
  {\it As a corollary, the system has monodromy if the holonomy of the heat connection
  $hol(\mathbb{\nabla}^{heat})$ is non trivial.}
  \end{enumerate}
     \end{theorem}
\begin{remarks}\par
\begin{enumerate} \item The BS-picture can be traded for the theta-picture:
the tracing of a non trivial path in $\pi_1(B)$ can be seen as a sort of adiabatic motion, causing
the variation of the basis of cycles and thence of the parameter $\tau$. The overall action on the
theta space is a sort of Berry phase (see e.g. \cite{CJ}), a signpost for monodromy.
The point is that in the theta-picture we are essentially acting in a single quantum Hilbert space
(in view of polarization independence).
This peculiarity pertains to the 2d-environment only.
We also stress the fact that the monodromy map $Q(Z_0)$  yields a {\it unitary operator}
(the crucial fact is that $\Im \tau$ does not change): this explains the notation $Q$.
If we read $Z_0$ classically, then  $Q(Z_0)$ is precisely its quantum counterpart and takes the form of a ``phase gate",
familiar from quantum computing in the qubit space ${\C{}}^2$ (\cite{CJ}).
We notice in passing that the appearance of a finite group like $\Z{}_4$ (cf. $Q(Z_0)^4 = Id_2$)  is to be expected on general grounds (cf. \cite{Loi,GH}).\par
\item We point out an important difference between our approach and Tyurin's one (\cite{tyurin1}):
in the latter case the BS-torus becomes the {\it real} part of an abelian variety;
in our case we have a 2d-BS-torus endowed with a complex structure. The latter is then holomorphically
quantized via $2$-level  theta functions, the natural substitute
for the BS-covariantly constant section whereupon the map  $\tau \mapsto \tau + 1$ acts
\`a la Berry. Hence there is no need of complexifying the manifold, study the ensuing complex monodromy
and then coming back to the (mechanically relevant) real picture (see also \cite{audin1, vivolo}).\par

\end{enumerate} 
\end{remarks} 

So, to summarize, the monodromy can be ascertained  via BS-wave functions, via
${\cal G}_0$-(in)equivalence and, in the 2d-case,
via the theta function description as well,
by means of the $SL(2, \Z{})$-action on the vector bundle determined by theta functions of level 2.\par

\subsection{Braiding and monodromy via theta functions}

In this Section we discuss the relationship between monodromy and the braiding of the Weiestra\ss \ roots.
More details on the braid group and its relationship with the modular group $PSL(2,\Z{})$ and, in particular, on related representations can be found e.g. in \cite{Mos, Tu, TW, BenveSpe}. \par
 
A  (faithful) representation of the braid group $B_3$  on  ${\C{}}^2$ via $SL(2,\Z{})$ can be realized by the matrices\par
\begin{equation}
b_1 = \begin{pmatrix} 
         1 & 0  \cr
         -1 & \hfill 1  \cr
        \end{pmatrix}        
\qquad 
b_2 = \begin{pmatrix} 
         1 & 1  \cr
         0 & 1  \cr
        \end{pmatrix}  
\end{equation}
Indeed one immediately checks the defining relation  $b_1\,b_2\,b_1 = b_2\,b_1\,b_2$. 
Notice that 
\begin{equation}
b_2 = (b_1)^{-T}
\end{equation}
\par
Also, one recognizes that
the trasformation $ \tau \mapsto \tau + 1$  can be represented by $b_2$. Thus the braid group generators are dual to each other,
from the point of view of classical-quantum monodromy (see Subsection 3).

This, in turn, can be read
on  the fundamental cycles, on the thetas, and
 on the Weierstra\ss \ roots of the associated elliptic curve:
$e_1 \mapsto e_1, e_2 \mapsto e_3, e_3 \mapsto e_2$ (see Subsection 2.3 and e.g. \cite{Tricomi}, 
 \cite{Mos}).
The following reference formulae are helpful in making this point:\par
\begin{equation}
2 \omega = \oint_{\gamma_1} {dz \over {\sqrt{P(z)} }} \qquad 
2 \omega^{\prime} = \oint_{\gamma_3} {dz \over {\sqrt{P(z)} }}
\end{equation}
with  $\omega^{\prime} = \tau \omega$, and where
 $\gamma_1$ surrounds $e_2$ and  $e_3$ - passing to the other sheet
of the ramified double cover, through the cut joining $e_3$ to $\infty$ -
and $\gamma_3$ encircles $e_1$ and  $e_2$ (cf. \cite{Tricomi}, fig. 19, p.85). 
In our context, Tricomi's $\gamma_3$ is our $\gamma_2$ (the changing cycle) whilst Tricomi's $\gamma_1$  is our $\gamma_1$ (the fixed cycle).
\par
The above considerations immediately lead us to the following

\begin{theorem}
{\it In the case of an isolated  focus-focus singularity, the variation of the rotation number is tantamount to the (multiple) switching of the roots
$e_2$ and $e_3$ (with the above conventions). More precisely,
if $[\gamma]$ is a generator for $\pi_1(B) \cong \Z{}$, then {\rm classical} monodromy is represented via}
\begin{equation}
m \cdot [\gamma] \leftrightarrow  {b_2}^{m}
\end{equation}
{\it whereas {\rm quantum} monodromy is given by
\begin{equation}
m \cdot [\gamma] \leftrightarrow  {b_1}^{m}
\end{equation}
In terms of the rotation number one obviously has}
\begin{equation}
\Theta \mapsto \Theta - {m}
\end{equation}
\end{theorem}

\section{The spherical pendulum revisited}
In this section we quickly point out a derivation of the monodromy of the spherical pendulum 
(the prototype of monodromic behaviour, see also e.g. \cite{bates}) via root braiding. We
refer to \cite{cushman-bates} for background and notation. 
The central object is the polynomial
\begin{equation}
P(x) = 2 (h - x) (1 - x^2) - j^2.
\end{equation}
The point $(j,h) = (0,1)$ is the only critical point in the  (punctured) open ``shield" (i.e. the $B$, in the present example).
Consider the circuit 
\begin{equation}
j =  \varepsilon \cos t , \quad h =  1 + \varepsilon \sin t, 
\end{equation}
for $t \in [0, 2\pi )$ and $ \varepsilon > 0$ small enough.
The roots of $P$ can be guessed via an $\varepsilon$-power series expansion, which immediately leads to the
 (exact!) expressions below:

\begin{equation}
x^{-} = -1 + \varepsilon^2 \,{{\cos^2 t} \over 8}, \quad x^{+} = 1 - {\varepsilon \over 2}\, (1 -\sin t),
\quad x^{0} = 1 + {\varepsilon \over 2}\, (1 +\sin t)
\end{equation}

As for $\tau = -\Theta + i T$, one has, for the spherical pendulum
\begin{equation}
\Theta = 2j \int_{x^{-}}^{x^{+}} {dx \over {(1 - x^2) \sqrt{P(x)} }}, \quad 
T = 2 \int_{x^{-}}^{x^{+}} {dx \over {\sqrt{P(x)} }}
\end{equation}

The first integral is a sum of elliptic integrals of the third kind, whereas the second one is of the first kind.
Now the basic result is
\begin{proposition}
 {\it $T$ is single-valued, whereas the variation of $\Theta$ along the above circuit equals $-1$.}
\end{proposition}
\noindent
This recovers monodromy for the spherical pendulum (\cite{duistermaat,cushman-bates}.\par
\smallskip
{\bf Sketch of Proof.} One has to study the braiding of the roots of $P$ along the closed path above;
now, the only delicate point is that at some positions the roots $x^{-}$, $x^{+}$
reach the limiting positions $-1$ and $+1$, and, in the first integral, one has to cope with
the onset of branching points, causing a rearrangement of the Riemann surface involved;
but this occurs in a non symmetrical way and induces the overall variation asserted above.
Nothing happens in the other case.\par

\section{Conclusions and outlook}
In this paper, acting within the framework provided by geometric quantization, 
we elaborated on the BS-wave function description of
 monodromy, by pointing out  its connection to gauge theory,
 and
via a theta function description as well,
upon viewing the $SL(2, \Z{})$-action on the vector bundle built up from theta functions of level 2.
We related monodromy to the braiding of the Weierstra\ss \ roots of the elliptic curve pertaining to the complex structure
$\tau$, again via theta functions. As a related application, we studied the monodromy of the spherical pendulum via braiding of the roots associated to
the elliptic integrals appearing in the expression of the relevant $\tau$. This method is, in principle, applicable to general completely integrable Hamiltonian systems
with two degrees of freedom.\par
Finally, we  expect  that a geometric quantization approach could be fruitful for dealing with the Chern-Duistermaat class as well.\par
\medskip

{\bf Acknowledgements.} The authors are grateful to A. Giacobbe and E. Previato for useful discussions.
The research of N.S. has been supported by a grant (AdR 819/07 ``Geometria globale dei sistemi completamente integrabili") 
from the Universit\`a degli
Studi di Verona, that of  M.S. by the Italian M.I.U.R. (ex 60\% funds).

\end{document}